\documentclass[12pt]{amsart}
\usepackage{amsmath,amsfonts}
\usepackage[mathscr]{eucal}
\usepackage{multirow}
\usepackage{amsthm,colordvi}
\usepackage{graphics}
\usepackage{graphicx}
\newcommand{\field}[1]{\ensuremath{\mathbb{#1}}}
\newcommand{\R}{\field{R}}

\renewcommand{\P}{\field{P}}
\newcommand{\Hess}{\mbox{\rm Hess}}

\newcommand{\QED}{\hfill$\Box$\medskip}

\newtheorem{teo}{Theorem}

\newtheorem{pro}[teo]{Proposition}

\title{Real Hessian Curves}

\author{Adriana Ortiz-Rodr\'\i guez}
\address{Instituto de Matem\'aticas,\\
     Universidad Nacional Aut\'onoma de M\'exico\\
     \'Area de la Inv. cient\'\i fica\\
      Circuito Exterior, Ciudad Universitaria\\
      M\'exico, D.F. 04510, M\'exico}
\email{aortiz@matem.unam.mx}
\author{Frank Sottile}
\address{Department of Mathematics\\
         Texas A\&M University\\
         College Station\\
         TX \ 77843\\
         USA}
\email{sottile@math.tamu.edu}
\urladdr{http://www.math.tamu.edu/\~{}sottile}
\thanks{Ortiz-Rodr\'\i guez supported in part by DGAPA and SEP-CONACyT project 41339-F}
\thanks{Sottile supported in part by NSF CAREER grant DMS-0538734}
\begin{document}
%%%%%%%%%%%%%%%%%%%%%%%%%%%%%%%%%%%%%%%%%%%%%%%%%%%%%%%%%

\begin{abstract}
 We give some real polynomials in two variables of degrees $4$, $5$, and
 $6$ whose hessian curves have more connected 
 components than had been known previously.
 In particular, we give a quartic polynomial whose hessian curve has 4 
 compact connected
 components (ovals), a quintic whose hessian curve has 8 ovals, and a sextic whose
 hessian curve has 11 ovals.
\end{abstract}

%%%%%%%%%%%%%%%%%%%%%%%%%%%%%%%%%%%%%%%%%%%%%%%%%%%%
\maketitle

\section*{Introduction}

The \Blue{\emph{parabolic curve}} on a generic smooth surface $S$
embedded in three-dimensional Euclidean space consists of the points 
where $S$ has zero Gaussian curvature. 
It separates \Blue{\emph{elliptic points}} (where the curvature
is positive) from \Blue{\emph{hyperbolic points}} (where the curvature is
negative).
These notions are well-defined for surfaces embedded in
affine or even projective space, as the sign of Gaussian curvature is
invariant under affine transformations.

If the surface $S$ is expressed locally as the graph $z=f(x,y)$ of a smooth
function $f$, then the sign of its hessian determinant
\[
   \Hess(f)\ :=\ 
    \left|\begin{matrix}f_{xx}&f_{xy}\\f_{yx}&f_{yy}\end{matrix}\right|
    \ =\ f_{xx}f_{yy}-f_{xy}^2\,,
\]
equals the sign of its curvature at the corresponding point.
Thus the parabolic curve is the image under $f$ of its \Blue{\emph{hessian curve}},
which is defined by $\Hess(f)=0$.  
When the surface $S$ is the graph of a polynomial $f\in\R[x,y]$,
this local description is global, and so questions about the disposition of the
parabolic curve on $S$ are equivalent to the same questions about the
hessian curve in $\R^2$.

Suppose that $d$ is even.
Harnack proved~\cite{ha} that a smooth plane curve of degree $d$ has
at most $1+\binom{d-1}{2}$ connected components in $\R\P^2$.
This is also the bound for the number of components of a compact curve in
$\R^2$ of degree $d$.
A non-compact curve in $\R^2$ of degree $d$ can have at most
$\binom{d-1}{2}$ bounded components (\Blue{{\it ovals}}) and $d$ unbounded
components.   
These unbounded components come from the intersection of the corresponding
curve in $\R\P^2$ with the line at infinity.
Harnack constructed a curve in $\R\P^2$ of degree $d$ with
$1+\binom{d-1}{2}$ components which has one component meeting the
line at infinity in $d$ points.
This Harnack curve shows that the bound for non-compact curves in $\R^2$ is
attained, and choosing a different line at infinity shows that the bound
for compact curves in $\R^2$ is also attained.

We are interested in the possible number and disposition of the components of the
hessian curve in $\R^2$ of a polynomial $f \in\R[x,y]$ of degree $n$.
This is problem 2001-1 in the list of Arnold's problems~\cite{ar1},
attributed to A.~Ortiz-Rodr\'\i guez.
See also the discussion of related problems 2000-1, 2000-2, 2001-1, and 2002-1.
The hessian of $f$ has degree at most $2n-4$.
By Harnack's Theorem, a compact hessian curve has at
most $(2n{-}5)(n{-}3){+}1$ ovals
and a non-compact hessian curve has at most $(2n{-}5)(n{-}3)$ ovals and
$2n{-}4$ unbounded components.

While we know of no additional restrictions on 
hessian curves, we are not assured that all possible configurations 
are acheived by hessians.
When $n$ is at least 4, simple parameter counting shows that not all
polynomials of degree $2n{-}4$ arise as hessians of polynomials of degree $n$. 
The placement of the set of hessian curves among all curves of degree 
$2n{-}4$ may restrict the possible configurations of hessian curves in
$\R^2$.
For example, a simple calculation shows that 
\[
   \Hess(f)\ =\ \left(\frac{(f_{xx}+f_{yy})}{2}\right)^2
  \  -\ \left(\frac{(f_{xx}-f_{yy})}{2}\right)^2\ -\ f_{xy}^2\,.
\]
Thus the hessian of a polynomial is a linear combination of 3 squares,
which shows that the  hessians lie in the second secant variety to the veronese
embedding of polynomials of degree $n{-}2$ in polynomials of degree
$n{-}4$ (the veronese consists of the perfect squares).

We also know of no general techniques for constructing hessian
curves with a prescribed configuration.
One of us (Ortiz-Rodr\'\i guez) investigated this
question~\cite{or1,or2} and constructed polynomials $f \in \R[x,y]$ of
degree $n$ whose hessians had $\binom{n-1}{2}$ ovals in $\R^2$. 
When $n$ is 4, 5, and 6, these numbers are 3, 6, and 10, respectively. 
We do not know if it is possible for a hessian curve to achieve the
Harnack bound, or more generally, which configurations are possible
for hessian curves. 

Here, we present a quartic polynomial $f$ whose hessian achieves the Harnack bound
of 4 ovals, a quintic whose hessian has 8 ovals, a sextic whose hessian has 11
ovals, as well as examples of non-compact hessian curves of quartics,
quintics, and sextics. 
These examples show that hessian curves can have more ovals than 
was previously known.
They were found in a computer search, using the software Maple.

Our method was to generate a random polynomial, compute its hessian, and
then compute an upper bound on its number of ovals in $\R\P^2$,
sometimes also screening for the number of unbounded components in $\R^2$.
This upper bound was one-half the minimum number of real critical points of a
projection to one of the axes, as each oval in $\R\P^2$ contributes at least
two critical points to the projection. 
We separately investigated compact hessian curves of sextics.
Polynomials whose upper bound for ovals was at least 4, 8, and 11 (for
quartics, quintics, and sextics, respectively) were saved for further
study.  
The further investigation largely involved viewing pictures in $\R^2$ of
these potentially interesting hessians. 
In all, only a few hundred polynomials warranted such further scrutiny.

We examined the hessians of 150 million quartics, 
40 million each of quintics and sextics, and over 200 million sextics with
compact hessians (the different protocol of pre-screening for compactness 
allowed a greater number to be examined).
This required 628 days of CPU time on several computers, most of which 
were running Linux on Intel Pentium processors with speeds
between 1.8 and 3 gigaHertz.
We did not find a quartic whose hessian had 3 ovals and 4 unbounded
components, 
nor a quintic whose hessian had more than 8 ovals in $\R\P^2$, nor a
sextic whose hessian had more than 11 ovals in $\R\P^2$. 
(The examples we give at the end with 12 ovals in $\R\P^2$ are pertubations
of a curve we found with 11 ovals.)
This suggests that it may not be possible for hessian curves in $\R^2$
to achieve the Harnack bounds.
Further pictures and computer code are at the web page\footnote{{\tt www.math.tamu.edu/\~{}sottile/stories/Hessian/index.html}.}.
%
%Quartic_Time:= 3471316.17:
%Quintic_Time:= 4267106.67:
%Sextic_Time :=32486915.45: 
%Compact_Sextic_Time:=14098316.83:
%  628 Days

 Tables~\ref{T1} and~\ref{T2} summarize this discussion concerning the
number of components of hessian curves.
The pairs $(o,u)$ in Table~\ref{T2} refer to ovals and unbounded
components, respectively.

%%%%%%%%%%%%%%%%%%%%%%%%%%%%%%%%%%%%%%%%%%%%%5
\begin{table}[htb]
\begin{tabular}{|l||c||c|c|c|c|c|}\hline
  Degree of $f$ & $n$  &   3 & 4 &5 & 6 &7\\\hline\hline
  Degree of hessian& $2n{-}4$ &2 & 4 & 6 &8  &10\\\hline 
  Harnack bound for hessian& $(2n{-}5)(n{-}3)+1$& 1 &4 &11&  22 & 37\\\hline
  Ortiz hessians~\cite{or1,or2}&$(n{-}1)(n{-}2)/2$& 1& 3 &6 &10& 15\\\hline
  New examples  &  & & 4 & 8 & 11 & --- \\\hline
\end{tabular}\vspace{5pt}
\caption{Ovals of compact hessian curves.}
 \label{T1}
\end{table}
%%%%%%%%%%%%%%%%%%%%%%%%%%%%%%%%%%%%%%%%%%%%%5

%%%%%%%%%%%%%%%%%%%%%%%%%%%%%%%%%%%%%%%%%%%%%5
\begin{table}[htb]
\begin{tabular}{|l||c||c|c|c|c|}\hline
  Degree of $f$ & $n$  &   3 & 4 &5 & 6 \\\hline\hline
  Degree of hessian& $2n{-}4$ &2 & 4 & 6 &8  \\\hline 
  Harnack bound&
  $((2n{-}5)(n{-}3),\,2n{-}4)$
  & (0,2) &(3,4) &(10,6)&  (21,8)\\\hline
\multirow{2}{90pt}{New examples}  & & &(2,4)  & (6,4) & (10,4) \\\cline{3-6}
  && & (3,2) & (7,2) & (11,2)   \\\hline
\end{tabular}\vspace{5pt}
 \caption{Configurations of non-compact hessians.}
 \label{T2}
\end{table}
%%%%%%%%%%%%%%%%%%%%%%%%%%%%%%%%%%%%%%%%%%%%%

\medskip 
{\bf Acknowledgments:} 
We want to thank to  V.I.~Arnold, L.~Ortiz-Bobadilla, V.~Kharlamov,
B.~Reznick, and E.~Ro\-sales-Gonz\'alez for their suggestions and comments.
Sottile thanks Universidad Nacional Aut\'onoma de M\'exico for its
hospitality, and where this investigation was initiated.

%%%%%%%%%%%%%%%%%%%%%%%%%%%%%%%%%%%%%%%%%%%%%%%%%%%%%%%%%%%%%%%%%%%%%%
\section*{Hessian curves with many ovals}

We begin with the following observation about hessian polynomials.

\begin{pro}
 A polynomial $h(x,y)$ is a hessian of some polynomial $f$ if and only if there exist 
 polynomials $p,q,r$ such that $p_y=q_x$,  $q_y=r_x$, and 
 $h=pr-q^2$.
\end{pro}

\noindent{\it Proof.}
 If $h$ is the hessian of $f$, then $h=f_{xx}f_{yy}-f_{xy}^2$, and $f_{xx}$, $f_{xy}$,
 and $f_{yy}$ satisfy these conditions.
 Conversely, if  $p$, $q$, and $r$ satisfy the conditions, then elementary
 integral calculus gives  
 polynomials $s$ and $t$ such that $s_x=p$,  $s_y=q$, $t_x=q$, and $t_y=r$.
 Since $s_y=t_x$, there is a polynomial $f$ with
 $f_x=s$ and $f_y=t$, and thus $h$ is the hessian of $f$.
\QED

\begin{teo}\label{t1}
There exists a real polynomial 
of degree $4$ in two variables whose hessian curve 
is smooth, compact, and consists of exactly four ovals.
\end{teo}

\begin{figure}[htb]  
\[
  \begin{picture}(280,153)
   \put(0,0){\includegraphics[width=10cm]{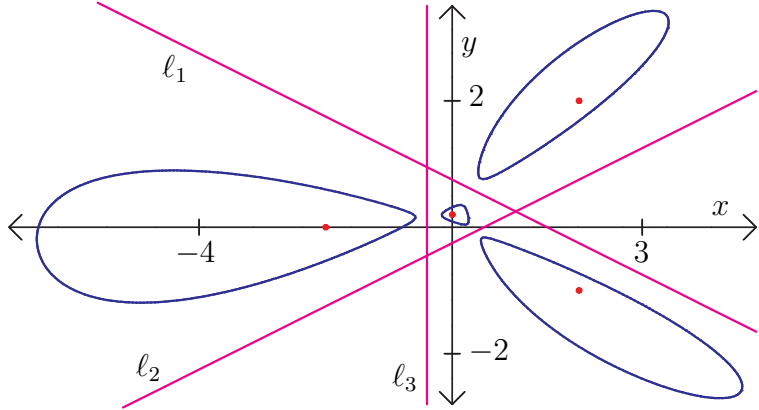}}
   \put(268,73){$x$}   \put(173,135){$y$}
   \put(65,55){$-4$} \put(238.5,55){$3$}
   \put(176,18){$-2$} \put(176,114){$2$}
   \put( 60,125){$\ell_1$}
   \put( 50,12){$\ell_2$}
   \put(147,08){$\ell_3$}
  \end{picture}
\]
    \caption{Quartic hessian curve with 4 compact components}
       \label{cuatro-ov}
\end{figure}

\noindent{\it Proof.}
Let $f$ be the polynomial
\[
   -2y^2 +2xy   +12x^2\ 
         +\ 10y^3  +3xy^2 -10x^2y  -13x^3\ 
       -\ 11y^4  +6xy^3  +9x^2y^2 -2x^3y -x^4\,.
\]
If we divide its hessian by $-4$, we obtain the polynomial
\begin{align*}
 h\ :=\  &\ \ 25-134x-374y+91x^2+948xy+1137y^2+429x^3+612x^2y-2313xy^2-876y^3\\
  &+63x^4+54x^3y-99x^2y^2-234xy^3+675y^4\,.
\end{align*}
We claim that the hessian curve, $h(x,y) = 0$, is a compact smooth curve
in $\R^2$ with exactly 4 connected components.
A visual proof is provided by the picture of the hessian curve in
Figure~\ref{cuatro-ov}.
%%%%%%%%%%%%%%%%%%%%%%%%%%%%%%%%%%%%%%%%%%%%%%%%%%%%%%%%%%%%%%%%%%%%%%%%%%%%%%%%
%%%%%%%%%%%%%%%%%%%%%%%%%%%%%%%%%%%%%%%%%%%%%%%%%%%%%%%%%%%%%%%%%%%%%%%%%%%%%%%%
This was drawn by Maple using its {\tt implicitplot} function with 
$200\times 200$ grid.
While this is sufficient to verify our claim, we provide alternative 
ad hoc  arguments.

We compute the values of the hessian at the four points inside
each oval of Figure~\ref{cuatro-ov}
\[
  h(-2,0)=-7068,\quad h(0,{\textstyle \frac{1}{5}})=-{\textstyle \frac{5124}{125}},\quad
  h(2,2)=-8508,\quad\mbox{and}\quad h(2,-1)=-6828\,.
\]
Next, we shall prove that $h$ is positive
on the line $\ell_\infty$ at infinity and on 
three lines of Figure~\ref{cuatro-ov},  
\[
  \ell_1\ \colon\  y=\frac{3}{4}-\frac{x}{2}, \qquad
  \ell_2\ \colon\  y=\frac{x}{2}-\frac{1}{4}, \quad\mbox{ and }\quad
  \ell_3\ \colon\  x=-\frac{2}{5}\,.
\]

The lines $\ell_\infty$, $\ell_1$, $\ell_2$, and $\ell_3$
divide $\R\P^2$ into 7 components. 
Since $h$ is  positive on these lines but is negative at the
four points $(-2,0)$, $(0,{\textstyle \frac{1}{5}})$, 
 $(2,2)$, and $(2,-1)$, which lie in different regions, the hessian curve
$h=0$ is compact and has at least one 1-dimensional component in each region
surrounding one of the four points. 
Since 4 is the maximum number of one-dimensional connected components of a quartic,
and such quartics are smooth, we deduce that the hessian curve is
smooth, compact, and consists of exactly four ovals.

For the claim about positivity, note that $h$ contains the
monomial term $63x^4$, and so it is positive at one point on the line at
infinity.
We show that $h$ does not vanish on any of these four lines, which 
implies our claim about positivity.
For this, we invoke a classical characterization of when a univariate
quartic has no real zeroes.
References may be found, for example in~\cite[\S71]{Dickson}.

Given a univariate quartic polynomial of the form
\[
  z^4 + \Blue{4}\alpha z^3 + \beta z^2 + \gamma z + \delta\,,
\]
linear substitution of $(z-\alpha)$ for $z$ gives the reduced quartic
\[
  z^4 + az^2+bz+c\,,
\]
where $a=\beta-6\alpha^2$, $b=\gamma-2\alpha\beta+8\alpha^3$, and 
$d=\delta-\alpha\gamma +\alpha^2\beta -3\alpha^4$.
The discriminant of this reduced quartic is
\[
   \Delta\ :=\ -4a^3b^2-27b^4+16a^4c-128a^2c^2+144ab^2c+256c^3\,.
\]
This criterion also uses the polynomial
\[
  L\ :=\ 2a(a^2-4c)+9b^2\,.
\]
Then the quartic has no real zeroes if and only if
 \begin{equation}\label{quartic_criterion}
  \Delta>0\quad\mbox{and either}\quad a\geq 0\quad\mbox{or}\quad L\geq
  0\ .
 \end{equation}

Homogenizing $h$, restricting it to the line at infinity, substituting
$y=1$, and dividing by 9 gives the quartic
\[
  q_\infty\ :=\  7x^4+6x^3-11x^2-26x+75\,.
\]

Restricting $h$ to the lines $\ell_1$, $\ell_2$, and $\ell_3$ 
and clearing denominators gives
\begin{eqnarray*}
  q_1 &:=&  21168x^4-157632x^3+592264x^2-337648x+58387\,,\\
  q_2 &:=&  20016x^4+4608x^3+377320x^2-278112x+52707\,,\makebox{ and}\\
  q_3 &:=& 421875y^4-489000y^3+1278975y^2-411710y+42073\,.
\end{eqnarray*}
These satisfy the criterion~\eqref{quartic_criterion} to have no real
zeroes, as may be seen from Table~\ref{discr}, where we give the
values of $\Delta$, $L$, and $a$, for each of these polynomials.\QED

%%%%%%%%%%%%%%%%%%%%%%%%%%%%%%%%%%%%%%%%%%%%%%%%%%%%%%%%%%%%%%%%%%%%%%%%%%%%%%%5
\begin{table}[htb]
\begin{tabular}{|c||c|c|c|c|}\hline
  Polynomial &$\Delta$&$L$&$a$\\\hline
  $q_\infty$&$\frac{5025022208}{16807}$&$\frac{564896}{2401}$&$\frac{-181}{98}$\raisebox{-8pt}{\rule{0pt}{22pt}}\\\hline
 $q_1$&$\frac{105415059013155058653376}{198607342807439307}$&$\frac{3692894126604316}{340405734249}$&$\frac{931453}{129654}$\raisebox{-8pt}{\rule{0pt}{22pt}}\\\hline
 $q_2$&$\frac{34807374069358185363904}{141964610099247963}$&$\frac{4123100447100116}{272136458889}$&$\frac{6549023}{347778}$\raisebox{-8pt}{\rule{0pt}{22pt}}\\\hline
 $q_3$&$\frac{10042565821320692218681168}{855261504650115966796875}$&$\frac{1376823939540422}{40045166015625}$&$\frac{1066423}{421875}$\raisebox{-8pt}{\rule{0pt}{22pt}}\\\hline
\end{tabular}\vspace{5pt}
\caption{Values of $\Delta$, $L$, and $a$.}\label{discr}
\end{table}
%%%%%%%%%%%%%%%%%%%%%%%%%%%%%%%%%%%%%%%%%%%%%%%%%%%%%%%%%%%%%%%%%%%%%%%%%%%%%%%
%
%   This is the output of a maple code.  The qudruple of number is 
%
%    [ Delta, a, L, 4c-a^2 ]
%
%7*x^4+6*x^3-11*x^2-26*x+75, 9
%       5025022208  -181  564896  101600
%      [----------, ----, ------, ------]
%         16807      98    2401    2401
%
%21168*x^4-157632*x^3+592264*x^2-337648*x+58387, 1/256
%      105415059013155058653376  931453  3692894126604316  39401665048
%     [------------------------, ------, ----------------, -----------]
%         198607342807439307     129654    340405734249     466948881
%
%20016*x^4+4608*x^3+377320*x^2-278112*x+52707, 1/256
%    34807374069358185363904  6549023  4123100447100116  -10299666301264
%   [-----------------------, -------, ----------------, ---------------]
%      141964610099247963     347778     272136458889      30237384321
%
%421875*y^4-489000*y^3+1278975*y^2-411710*y+42073, 1/625
%   10042565821320692218681168  1066423  1376823939540422  -3304239016591
%  [--------------------------, -------, ----------------, --------------]
%    855261504650115966796875   421875    40045166015625    533935546875
%
%%%%%%%%%%%%%%%%%%%%%%%%%%%%%%%%%%%%%%%%%%%%%%%%%%%%%%%%%%%%%%%%%%%%%%%%%%%%%%%5

Each of the remaining curves we discuss is smooth,\rule{0pt}{15pt} each oval has
exactly two vertical and two horizontal tangents, and each unbounded
component has exactly one vertical and one horizontal tangent.
These claims are best verified symbolically.
For each, we give the polynomial $f$ and display a picture of the
hessian curve, drawn with the {\tt implicitplot} function of Maple.
These were rendered, at least locally, with a grid size sufficiently
small to separate the tangents, and therefore provide a faithful picture
of the hessian curves as curves in $\R^2$.

\begin{figure}[htb] 
\[  
  \begin{picture}(315,135)(0,-20)
   \put(0,0){\includegraphics[width=4cm]{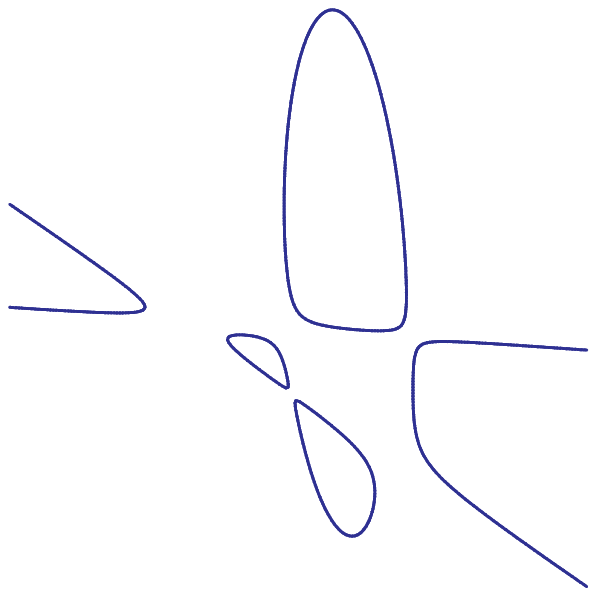}}
    \put(55,-20){(a)}
   \put(185,0){\includegraphics[width=4cm]{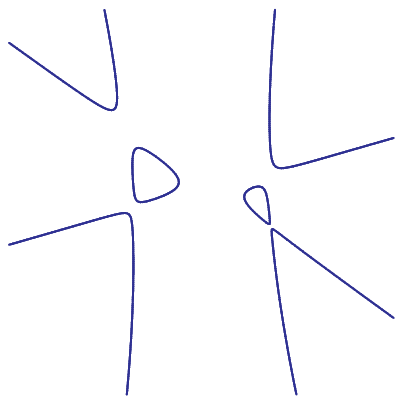}}
    \put(235,-20){(b)}
  \end{picture}
\]
 \caption{Hessians of quartics\label{F:two}}
\end{figure}

Figure~\ref{F:two}(a) displays the hessian curve of the quartic
\[
    22x^2+36xy+24y^2
  -80x^3-10x^2y+71xy^2 +39y^3
  +15x^4+4x^3y-3x^2y^2-21xy^3-17y^4\,,
\]
which has 3 ovals and 2 unbounded components.
Figure~\ref{F:two}(b) displays the hessian curve of the quartic
\[
   -70x^2 -35xy -2y^2 - 93x^3 -14x^2y +41xy^2-70y^3
   + 31x^4 +7x^3y -30x^2y^2 +37xy^3 +91y^4\,,
\]
which has 2 ovals and 4 unbounded components.
%%%%%%%%%%%%%%%%%%%%%%%%%%%%%%%%%%%%%%%%%%%%%%%%%%%%%%%%%%%%%%%%%%%%%%%%%%%%%%%%
%%%%%%%%%%%%%%%%%%%%%%%%%%%%%%%%%%%%%%%%%%%%%%%%%%%%%%%%%%%%%%%%%%%%%%%%%%%%%%%%
While we have generated and checked 150 million quartics, we did not
find one whose hessian curve achieves the Harnack bound of 3 ovals and
4 unbounded components. 

\begin{figure}[htb] 
\[   
  \begin{picture}(350,170)(0,-20)
   \put(0,0){\includegraphics[width=5cm]{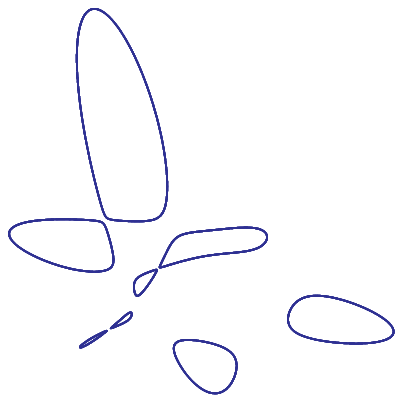}}
    \put(65,-20){(a)}
    \put(200,0){\includegraphics[width=5cm]{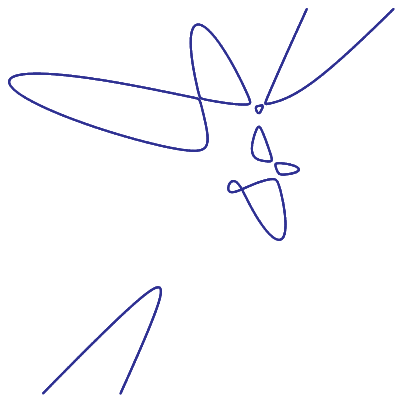}}
    \put(265,-20){(b)}
  \end{picture}
\]
 \caption{Hessians of quintics\label{F:three}}
\end{figure}

Figure~\ref{F:three}(a) displays the hessian curve of the quintic
\begin{align*}
   &4y^2+xy-6x^2\ -\ 25y^3+24xy^2+15x^2y-33x^3\ +\ y^4-3xy^3+15x^2y^2-19x^3y-26x^4\\
   & +\ 33y^5-2xy^4-23x^2y^3-30x^3y^2-26x^4y+31x^5\ ,
\end{align*}
which is compact with 8 ovals.
%%%%%%%%%%%%%%%%%%%%%%%%%%%%%%%%%%%%%%%%%%%%%%%%%%%%%%%%%%%%%%%%%%%%%%%%%%%%%%%%
%%%%%%%%%%%%%%%%%%%%%%%%%%%%%%%%%%%%%%%%%%%%%%%%%%%%%%%%%%%%%%%%%%%%%%%%%%%%%%%%

Figure~\ref{F:three}(b) displays the hessian curve of the quintic
\begin{align*}
&-54y^2-103xy-26x^2\ -\ 88y^3+45xy^2+91x^2y-96x^3 -\ 12y^4+43xy^3+6x^2y^2\\
&+\ 11x^3y+49x^4\ \ +\ \ 22y^5-20xy^4-38x^2y^3-14x^3y^2+45x^4y+76x^5\ ,
\end{align*}
which has 7 ovals and 2 unbounded components.

\begin{figure}[htb] 
\[
  \begin{picture}(310,140)
   \put(0,0){\includegraphics[width=2.6cm]{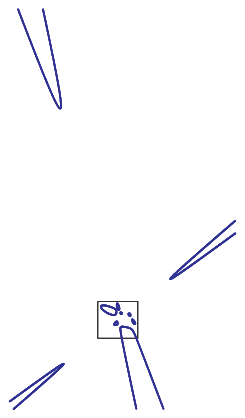}}
   \put(160,0){\includegraphics[width=5.195cm]{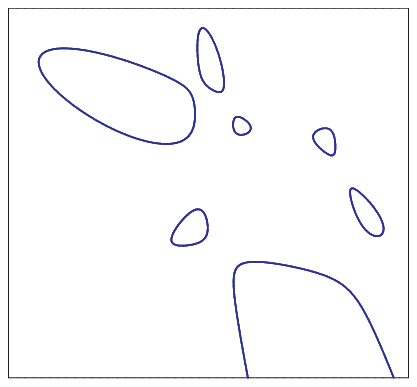}}
  \end{picture}
\]
 \caption{Hessian of a quintic with 6 ovals and 4 unbounded components\label{F:four}}
\end{figure}

Figure~\ref{F:four} displays the hessian curve of the quintic
\begin{align*}
 &60y^2 + 21xy   + 76x^2
+95y^3 - 18xy^2 - 79x^2y   + 88x^3
-25y^4 - 22xy^3 + 50x^2y^2 - 9x^3y    - 5x^4\\
&-57y^5 - 50xy^4 + 21x^2y^3 + 87x^3y^2 + 35x^4y - 56x^5\ ,
\end{align*}
which has 6 ovals and 4 unbounded components.
The boxed region on the left has been expanded in the picture on the right.
%%%%%%%%%%%%%%%%%%%%%%%%%%%%%%%%%%%%%%%%%%%%%%%%%%%%%%%%%%%%%%%%%%%%%%%%%%%%%%%%
%%%%%%%%%%%%%%%%%%%%%%%%%%%%%%%%%%%%%%%%%%%%%%%%%%%%%%%%%%%%%%%%%%%%%%%%%%%%%%%%

These quintics all have 8 ovals in $\R\P^2$.
While we have generated and checked 40 million quintics, we did not 
find any with more ovals.\medskip

\begin{figure}[htb] 
\[
  \begin{picture}(400,150)
   \put(0,0){\includegraphics[width=6cm]{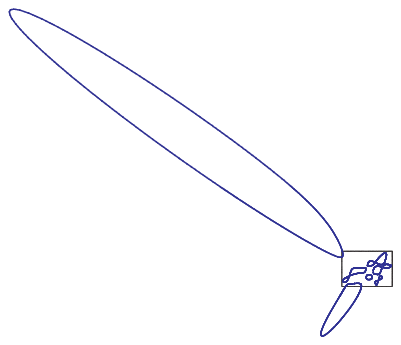}}
   \put(230, 20){\includegraphics[width=5.5cm]{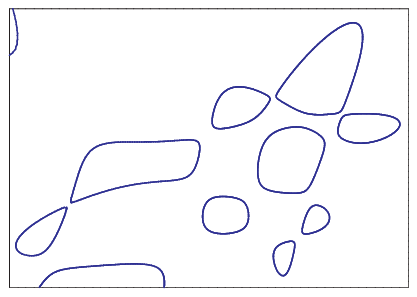}}
  \end{picture}
\]
 \caption{Hessian of a sextic with 11 ovals.\label{New}}
\end{figure}

%%%%%%%%%%%%%%%%%%%%%%%%%%%%%%%%%%%%%%%%%%%%%%%%%%%%%%%%%%%%%%%%%%%%%%%%%%%%%%%%
Figure~\ref{New} displays the hessian curve of the sextic
\begin{align*}
 &   45y^2-47xy-30x^2\ 
  +\ 96y^3-xy^2+8x^2y+54x^3\\
 & -\ 96y^4-64xy^3-50x^2y^2-33x^3y+91x^4\\
 & -\ 100y^5+84xy^4-43x^3y^2+66x^4y-58x^5\\
 & +\ 70y^6+90xy^5-28x^2y^4-53x^3y^3+43x^4y^2+36x^5y-38x^6\ ,
\end{align*}
which has 11 ovals.
The boxed region on the left has been expanded in the picture on the right.
%%%%%%%%%%%%%%%%%%%%%%%%%%%%%%%%%%%%%%%%%%%%%%%%%%%%%%%%%%%%%%%%%%%%%%%%%%%%%%%%
%%%%%%%%%%%%%%%%%%%%%%%%%%%%%%%%%%%%%%%%%%%%%%%%%%%%%%%%%%%%%%%%%%%%%%%%%%%%%%%%

\begin{figure}[htb] 
\[  
  \begin{picture}(435,225)(0,-20)
   \put(0,0){\includegraphics[width=7cm]{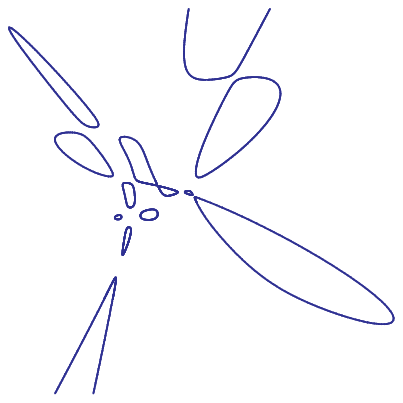}}
    \put(91,-20){(a)}
    \put(230,0){\includegraphics[width=7cm]{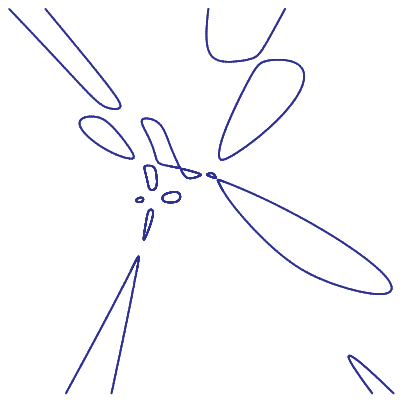}}
    \put(331,-20){(b)}
  \end{picture}
\]
 \caption{Hessians of sextics\label{F:five}}
\end{figure}

Figure~\ref{F:five}(a)
displays the hessian curve of the sextic
\begin{align*}
&-\ 53y^2-31xy+59x^2\ \ 
 -\ 79y^3+82xy^2-52x^2y+22x^3\\
&+\ 75y^4-27xy^3+63x^2y^2-85x^3y-89x^4\\
&+\ 80y^5+27xy^4-69x^2y^3+17x^3y^2-7x^4y-43x^5\\
&-\ 25y^6+17xy^5+27x^2y^4-55x^3y^3-37x^4y^2+59x^5y+45x^6\ ,
\end{align*}
which has 11 ovals and 2 unbounded components.
%%%%%%%%%%%%%%%%%%%%%%%%%%%%%%%%%%%%%%%%%%%%%%%%%%%%%%%%%%%%%%%%%%%%%%%%%%%%%%%%
%%%%%%%%%%%%%%%%%%%%%%%%%%%%%%%%%%%%%%%%%%%%%%%%%%%%%%%%%%%%%%%%%%%%%%%%%%%%%%%%

Figure~\ref{F:five}(b)
displays the hessian curve of the sextic
\begin{align*}
&-\ 80y^2-46xy+89x^2\ \ 
-\ 118y^3+123xy^2-78x^2y+33x^3\\
&+\ 113y^4-40xy^3+94x^2y^2-128x^3y-133x^4\\
&+\ 120y^5+40xy^4-104x^2y^3+25x^3y^2-10x^4y-64x^5\\
&-\ 37y^6+25xy^5+40x^2y^4-82x^3y^3-56x^4y^2+89x^5y+67x^6\ ,
\end{align*}
which has 10 ovals and 4 unbounded components.
Both hessian curves have 12 ovals in $\R\P^2$.

Despite examining over 240 million sextics, we did not find
any sextics whose hessian curves had more than 11 ovals in $\R\P^2$.
These last two examples, which have 12 ovals in $\R\P^2$, 
are perturbations of a sextic found
in the search whose hessian curve had 11 ovals in $\R\P^2$.

%%%%%%%%%%%%%%%%%%%%%%%%%%%%%%%%%%%%%%%%%%%%%%%%%%%%%%%%%%%%%

\end{document}